\documentclass[12pt]{article}
\usepackage{amsmath,amssymb,xcolor}
\usepackage[T1]{fontenc}

\newcommand{\eq}[1]{(\ref{eq:#1})}
\renewcommand{\d}{\text{d}}

\newtheorem{Theorem}{Theorem}
\newtheorem{Lemma}{Lemma}

\newtheorem{Remark}{Remark}

\def\cadlag{c\`adl\`ag\ }

\usepackage[a4paper,bindingoffset=0.5cm,left=2cm,right=2cm,top=2.5cm,bottom=2cm,footskip=.8cm]{geometry}

\begin{document}

\title{On the area between a L\'evy process with secondary jump inputs and its reflected version}

\author{
Offer Kella\thanks{Department of Statistics and Data Science, The Hebrew University of Jerusalem; Jerusalem 9190501, Israel. {\tt offer.kella@gmail.com}} 
\thanks{Partially supported by the Vigevani Chair in Statistics.} \ and Michel Mandjes\thanks{Mathematical Institute, Leiden University, P.O. Box 9512,
2300 RA Leiden,
The Netherlands. {\tt m.r.h.mandjes@math.leidenuniv.nl}} \thanks{Partially supported by the NWO Gravitation project N{\sc etworks}, grant number 024.002.003.}
}

\date{\today}

\maketitle

\begin{abstract}
\noindent
We study the stochastic properties of the area under some function of the difference between (i) a spectrally positive L\'evy process $W_t^x$ that jumps to a level $x>0$ whenever it hits zero, and (ii) its reflected version $W_t$. Remarkably, even though the analysis of each of these areas is challenging, we succeed in attaining explicit expressions for their difference. The main result concerns the Laplace-Stieltjes transform of the integral $A_x$ of (a function of) the distance between $W_t^x$ and $W_t$ until $W_t^x$ hits zero. This result is extended in a number of directions, including the area between $A_x$ and $A_y$ and a Gaussian limit theorem. We conclude the paper with an inventory problem for which our results are particularly useful.
\end{abstract}


\bigskip
\noindent {\bf AMS Subject Classification (MSC2020):} Primary 60G51; Secondary 60K25.

\section{Introduction}
Arguably one of the most intensively studied objects in stochastic operations research is the conventional M/G/1 queue: customers arrive according to a Poisson process, bring along i.i.d.\ amounts of work, which is put into a buffer that is emptied at a unit rate. 
It has been realized that many results for the workload of an M/G/1 system carry over to that of a storage system, in this paper denoted by $W\equiv\{W_t\}_{t\geqslant 0}$, driven by a spectrally-positive L\'evy process \cite{kyprianou} (i.e., a L\'evy process with no negative jumps). 
One of the key results that has been obtained for this class of reflected L\'evy processes is the {\it Pollaczek-Khinchine formula}, providing the Laplace-Stieltjes transform of the storage level in stationarity (as well as its counterpart after an exponentially distributed time); see e.g.\ \cite{DM}.  

In the M/G/1 system, whenever the system becomes empty, the server begins an idle time which is exponentially distributed after which another busy period begins. Alternatively, we could let the server start performing some secondary work during times at which it is idle. This could be preventive maintenance, or, more generally, any type of low-priority work. In the model considered in the present paper, this additional work, any time the buffer drains, has size $x>0$. Note that such a mechanism also leads to a meaningful model when the driving process is a spectrally-positive L\'evy process, in the sequel denoted by $W^x\equiv\{W_t^x\}_{t\geqslant 0}$.

Under a cost specific structure, which will be defined in this paper, it is relevant to decide whether such a secondary work is advisable, and, if it is, what the optimal value of the `pulse' of secondary work $x$ is.
To this end, with $T_x$ denoting the first time the storage level hits $0$ given it started at $x>0$, the object that needs to be studied is
\begin{equation}
A_x:=\int_0^{T_x}h(W^x_t-W_t)\,\d t
\end{equation}
for $x>0$ and a non-negative Borel function $h(\cdot)$. 
Areas such as $\int_0^{T_x}W_t\,\d t$ and $\int_0^{T_x}W^x_t\,\d t$ are notoriously hard to analyze, even in the M/G/1 case, which makes it rather remarkable that we succeed in finding a clean and explicit expression for the Laplace-Stieltjes transform of $A_x$.

The workload model with secondary input has been described in \cite{kw92} and, for the special case of Brownian motion, in \cite{kw90}. More general L\'evy processes with secondary jump inputs have been considered in \cite{kw91,bk14}. Recently, results for the moments of $\int_0^{T_x}W^x_t\,\d t$ have been derived in \cite{GJM}, but only for a subclass of the spectrally-positive L\'evy processes. 
In the setup of a G/G/1 queue, the first two moments of the integral of the storage process during a busy period are given in \cite{Borovkov2003,Iglehart1971}. For the M/M/1 queue, \cite{Borovkov2003} identified the Laplace-Stieltjes transform of this area in terms of Whittaker functions.  Some excellent textbooks on L\'evy processes are \cite{applebaum,bertoin,kyprianou,sato}. 

This paper is organized as follows. In Section~\ref{S2} we include a brief account of spectrally-positive L\'evy processes and associated hitting times, and provide a precise definition of all relevant objects. Then the main result, providing the Laplace-Stieltjes transform of $A_x$, is given in Section~\ref{S3}, including a series of ramifications. A procedure to recursively determine higher moments can be found in Section~\ref{S4} and a functional Gaussian limit theorem in Section~\ref{S5}. Finally, in Section~\ref{S6} we study an inventory problem with a specific underlying cost structure for which our results are particularly useful.

\section{Setup and preliminaries}\label{S2}
In what follows, let {\em LST, LS}, and {\em a.s.}~abbreviate {\em Laplace-Stieltjes transform, Laplace-Stieltjes}, and {\em almost surely}, respectively.

\subsection{Spectrally-positive L\'evy processes.}\label{SPL}
Let $X\equiv\{X_t\}_{t\geqslant 0}$ be a spectrally-positive L\'evy process starting from zero. 
By this we mean that $X$ is a \cadlag process having stationary and independent increments without any negative jumps, with $X_0=0$. 
This type of L\'evy processes is uniquely determined by a triplet $(c,\sigma^2,\nu)$, where $c\in\mathbb{R}$ and $\sigma^2\geqslant 0$ are constants, and $\nu(\cdot)$ is a measure on the Borel sets of $\mathbb{R}$, usually referred to as the {\em L\'evy measure}, which satisfies 
\begin{equation}\nu(-\infty,0]=0,\:\:\:\nu(1,\infty)<\infty,\:\:\mbox{and}\:\:\int_{(0,1)}x^2\,\nu(\d x)<\infty.\end{equation}
In this spectrally-positive case it is well known that ${\mathbb E}\,e^{-\alpha X_t}<\infty$ for every $\alpha,t\geqslant 0$ and moreover that ${\mathbb E}\,e^{-\alpha X_t}=e^{\varphi(\alpha)t}$, where $\varphi:[0,\infty)\to \mathbb{R}$ is called the LS exponent of $X$ and has the form
\begin{equation}
\varphi(\alpha)=-c\alpha+\frac{\sigma^2\alpha^2}{2}+\int_{(0,\infty)}\left(e^{-\alpha x}-1+\alpha x\,1_{(0,1]}(x)\right)\nu(\d x)\,.
\end{equation}
It is also known that $\varphi(\cdot)$ is a convex function, and unless $X$ is nondecreasing (in which case $X$ is called a {\em subordinator}) $\varphi(\alpha)$ converges to infinity as $\alpha\to\infty$. 
In fact, $\varphi(\cdot)$ is differentiable of any order on $(0,\infty)$, and $\varphi^{(n)}(0)\equiv \lim_{\alpha\downarrow 0}\varphi^{(n)}(\alpha)$ has the following form 
\begin{equation}
\varphi^{(n)}(0)=
\begin{cases}
-c-\int_{(1,\infty)}x\,\nu(\d x),&n=1,\\
\sigma^2+\int_{(0,\infty)}x^2\,\nu(\d x),&n=2,\\
(-1)^n\int_{(0,\infty)}x^n\,\nu(\d x),&n\in\{3,4,\ldots\}\,,
\end{cases}
\end{equation}
where it is noted, for any $n$, that $\varphi^{(n)}(0)$ is not necessarily finite. 
In particular we have that ${\mathbb E}\,X_t=-\varphi'(0)\,t$. In addition, if $\varphi'(0)$ is finite, then $\text{Var}\,X_t=\varphi''(0) \,t$.

If $\int_{(0,1]}x\,\nu(\d x)<\infty$, then the L\'evy process $X$ can be decomposed as $X_t=J_t+B_t$, where $J$ and $B$ are independent processes, $J$ is a pure jump subordinator with L\'evy measure $\nu(\cdot)$ and $B$ is a Brownian motion with drift $c-\int_{(0,1]}x\,\nu(\d x)$ and variance coefficient $\sigma^2\geqslant 0$. 
The only case for which $X$ is a subordinator is when $\sigma^2=0$ and $c\geqslant \int_{(0,1]}x\,\nu(\d x)$. 
This is also the only case for which $\varphi(\alpha)\leqslant 0$ for all $\alpha\geqslant 0$.
Throughout this paper we assume that in addition to $X$ being a \cadlag spectrally-positive L\'evy process, it is {\em not a subordinator}. 
We equivalently have that $\varphi(\alpha)\to\infty$ as $\alpha\to\infty$, which is also equivalent to the assumption that either $\sigma^2>0$ or $c<\int_{(0,1]}x\,\nu(\d x)$, which includes the case $\int_{(0,1]}x\,\nu(\d x)=\infty$.

It should be mentioned that the only case where $\varphi(\cdot)$ is not {\it strictly} convex is when $\sigma^2=0$ and $\nu(0,\infty)=0$. 
In this case $X$ corresponds to a linear drift: we a.s.\ have that $X_t=ct$ with $c<0$ (since $X$ is not a subordinator), and $\varphi(\alpha)=-c\alpha$. 
In all other cases $\varphi(\cdot)$ is strictly convex, hence, eventually strictly increasing: (i)~if $\varphi'(0)\geqslant 0$, then $\varphi(\cdot)$ is strictly increasing on $[0,\infty)$, (ii)~if $\varphi'(0)<0$ (including the case $\varphi'(0)=-\infty$), then $\varphi(\cdot)$ first decreases (thus attaining negative values), and then (strictly) increases to infinity. 
This means that if we denote 
\begin{equation}\alpha_0:=\inf\{\alpha\,|\,\varphi(\alpha)>0\}\end{equation}
(which is zero when $\varphi'(0)\geqslant 0$ and strictly positive otherwise), then $\varphi:[\alpha_0,\infty)\mapsto [0,\infty)$ is a strictly increasing and continuous on its domain, and therefore has an inverse $\varphi^{-1}:[0,\infty)\mapsto [\alpha_0,\infty)$. 

\subsection{Hitting times, local time}
For $x\geqslant 0$, consider the hitting time
\begin{equation}
T_x:=\inf\{t\geqslant 0\,|\,X_t<-x\}\,.
\end{equation}
It is well known (e.g., \cite[Th.~3.12, Cor.~3.13, Cor.~3.14]{kyprianou}) that the process $T\equiv \{T_x\}_{x\geqslant 0}$ is (right continuous) subordinator with LS exponent $-\varphi^{-1}(\cdot)$.
When $\varphi'(0)<0$, then $\varphi^{-1}(0)=\alpha_0>0$ and ${\mathbb P}(T_x<\infty)=e^{-\alpha_0 x}$, so that in this case the process $T$ is a {\it killed subordinator} with killing rate $\alpha_0$.
If $\varphi'(0)=0$, then we have that ${\mathbb P}(T_x<\infty)=1$ but ${\mathbb E}\,T_x=\infty$ (for any $x\geqslant 0$). 
If $\varphi'(0)>0$, then ${\mathbb E}\,T_x=x/\varphi'(0)<\infty$.
In this case all the finite moments of $T_x$ are expressions involving the derivatives of $\varphi^{-1}(\alpha)$ at zero, which in turn are expressions involving $\varphi^{(n)}(0)$ for $n\in\{ 0,1,\ldots\}$ via an application of Fa\`{a} di Bruno's formula. 

Now, consider the reflected process $W\equiv\{W_t\}_{t\geqslant 0}$ associated with $X$. It can be written as the sum of the {\it net input process} $X\equiv\{X_t\}_{t\geqslant 0}$ and the {\it local time} $L\equiv\{L_t\}_{t\geqslant 0}$:
\begin{equation}
W_t=X_t+L_t\,.
\end{equation}
where $L_t=-\inf_{0\leqslant s\leqslant t}X_s$ is the unique right continuous process with $L_0=0$, $L$ is nondecreasing, $W$ is non-negative  and \cite{k06}
\begin{equation}
\{t\,|\,W_t=0\}\subset \{t\,|\,\forall s\in[0,t): L_s<L_t\}\,.
\end{equation}
For the case that $X$ is a compound Poisson with unit negative drift, the local time $L$ can be interpreted as the cumulative idle time of the corresponding M/G/1 queue.
It increases linearly with rate one during idle times and is constant during busy times. When, a bit more generaly, $X$ is the difference between some subordinator and a positive linear drift of which rate is often referred to as the server's or station's {\em capacity}, $L$ can be interpreted as the {\em cumulative unused capacity}. In all other cases belonging to the class of L\'evy processes introduced in Section \ref{SPL}, $L$ is a more complicated process. 

Since there are no negative jumps, it is easy to check that $L$ is continuous. In addition, there are the following three equivalent characterizations involving the processes $L$ and $T$.

\begin{Lemma}\label{L1}
The following three assertions are equivalent, for any $x,t\geqslant 0$:
\begin{enumerate}
\item $T_x\geqslant t$.
\item $X_s\geqslant -x$ for all $0\leqslant s\leqslant t$.
\item $L_t\leqslant x$.
\end{enumerate}
\end{Lemma}

\noindent {\it Proof.}
First note that $2\Leftrightarrow 3$ follows by definition. 
If $X_s\geqslant -x$ for all $0\leqslant s\leqslant t$, then $T_x$ cannot be smaller than $t$ and hence $\{L_t\leqslant x\}\subset \{T_x\geqslant t\}$, so $2\Rightarrow 1$. 
Now, if $L_t>x$, then by the continuity of $L$ it follows that $L_{t-}>x$ and thus
\begin{equation}
T_x=\sup\{s\,|\,\forall u\in[0,s]:X_u\geqslant -x\}=\sup\{s\,|\,L_s\leqslant x\}<t 
\end{equation}
Therefore, $\{L_t>x\}\subset \{T_x<t\}$, which implies that $\{L_t\leqslant x\}\supset \{T_x\geqslant t\}$, so $1\Rightarrow 2$. 
As an aside, we note that when $X$ is not spectrally positive, then $L$ is not (a.s.) continuous and the inclusion $\{L_t>x\}\subset \{T_x<t\}$ can fail. \hfill$\Box$

\medskip

The equivalence $1\Leftrightarrow 3$ in Lemma \ref{L1} means that the local time process $L$ is the generalized inverse of the hitting time process $T$.
More concretely, on one hand we have that
\begin{equation}
L_t=\inf\{x\,|\,L_t\leqslant x\}=\inf\{x\,|\,T_x\geqslant t\}=\sup\{x\,|\,T_x<t\}\,,
\end{equation}
and on the other hand
\begin{equation}
T_x=\sup\{t\,|\,T_x\geqslant t\}=\sup\{t\,|\,L_t\leqslant x\}=\inf\{t\,|\,L_t>x\}\,.    
\end{equation}
Also, it it is easy to check that $L_{T_x}=-X_{T_x}=x$ (hence, also $W_{T_x}=0$), but in general $T_{L_t}\geqslant t$ and equality holds only for $t$ such that $L_s<L_t$ for all $s<t$. 

\medskip

Let us recall that if $F$ is a cumulative distribution function of some random variable $X$ and, for $u\in(0,1)$, $G(u)=\inf\{x\,|\,F(x)\geqslant u\}$ is its generalized inverse, then for $U\sim\text{Uniform}(0,1)$ we have that $G(U)$ is distributed according to $F$. 
Therefore for any non-negative  Borel function $f(\cdot)$ we have that
\begin{equation}
    \int_{\mathbb{R}}f(x)F(\d x)={\mathbb E}\,f(X)={\mathbb E}\,f(G(U))=\int_0^1 f(G(u))\d u\,.
\end{equation}
If $0<T_x<\infty$, taking $F(y)=T_y/T_x$ for $y\in[0,x]$, then it directly follows that $G(u)=L_{T_xu}$ and thus
\begin{equation}
\int_0^{T_x}f(L_t)\,\d t=T_x\int_0^1f(\underbrace{L_{T_xu}}_{G(u)})\,\d u=T_x\int_0^xf(y)\underbrace{T_{\d y}/T_x}_{F(\d y)}=\int_0^xf(y)\,T_{\d y}\,.
\end{equation}
If $T_x=0$, then $T_y=0$ for all $y\in[0,x]$ and hence both sides are zero. 
Therefore, if $\varphi'(0)\geqslant 0$ (so that ${\mathbb P}(T_x<\infty)=1$), then we have, for every $x\geqslant 0$ and every non-negative  Borel function, that
\begin{equation}\label{eq:LT}
\int_0^{T_x}f(L_t)\,\d t=\int_0^xf(y)\,T_{\d y}\,.
\end{equation}

As can be found in e.g.\ \cite{kw92}, if $\varphi'(0)>0$, then the process $W$ has a stationary distribution (which equals the process' limiting distribution, which equals the process' ergodic distribution), which from here on will be called a {\em steady-state distribution}.
It has the LST $\alpha\varphi'(0)/\varphi(\alpha)$ for $\alpha\geqslant 0$, which is sometimes referred to as the {\it generalized Pollaczeck-Khinchine formula} \cite{DM}. 
Therefore from here on we will assume that $\varphi'(0)>0$. Under this assumption it is impossible that $X$ is a subordinator and thus we do not need to explicitly assume anything else.

We conclude this section by considering, for a given $x>0$, the process $S^x\equiv \{S_n^x\}_{n\in{\mathbb N}}$ defined through $S^x_0=0$ and, for $n\ge 1$,
\begin{equation}S_n^x:=\inf\{t\,|\,X_t<-nx\}.\end{equation} This is a random walk with independent and identically distributed non-negative increments, each of them distributed like $T_x$. Also let $W^x\equiv \{W_t^x\}_{t\geqslant 0}$ be defined by
\begin{equation}
W^x_t=X_t+x\sum_{n=1}^\infty n\,1_{\{[S^x_{n-1},S^x_n)\}}(t)\,. 
\end{equation}
The main idea behind introducing this process is that it models `secondary jump inputs': the process $W^x$ is constructed such that whenever it is about to cross zero (from above), then it immediately jumps to the level $x>0$. 
We clearly have that $W^x_t\geqslant W_t$ for every $t\geqslant 0$ provided we use the same driving process $X$ in both, and it is known \cite{kw92} that the steady-state distribution of $W^x$ is a convolution or the steady-state distribution of $W$ and a $\text{Uniform}(0,x)$ distribution.

\begin{Remark}\rm
    We note in passing that when $\varphi'(0)>0$ we have that, for $x>0$, $\lim_{t\to\infty}L_t/t=\varphi'(0)=x/{\mathbb E}T_x$. As was mentioned earlier, for the case where the L\'evy process is a subordinator minus some positive drift $d$, $L_t$ is the cumulative unused capacity until time $t$ and therefore $\varphi'(0)$ is the amount of unused capacity per unit time. $\varphi'(0)=x/{\mathbb E}T_x$ is also the amount of secondary work which is performed per unit time. This makes intuitive sense as for the process $W^x_t$ all unused capacity is used for secondary work.\hfill$\Diamond$
\end{Remark}

\section{Distribution of an area}\label{S3}
In this section we explore the distribution of the area underneath a function of the secondary work content, that is, the difference between a spectrally-positive L\'evy process with secondary jump inputs $W^x$ and its reflected version $W$.
More precisely, bearing in mind that $T_x$ is a regeneration epoch for both $W^x$ and $W$, we set out to investigate the characteristics of the random variable
\begin{equation}
A_x:=\int_0^{T_x}h(W^x_t-W_t)\,\d t
\end{equation}
for $x>0$ and a non-negative Borel function $h(\cdot)$. We recall (e.g., \cite[Eqn.~(8)]{bkm19}, \cite[Th.~6.2.5]{lukacs}, \cite[Lem.~1.1]{jv83}), that for a subordinator $J$ with LS exponent $-\eta(\cdot)$, we have that for any non-negative Borel function $f(\cdot)$
\begin{equation}\label{eq:eta}
{\mathbb E}\,\exp\left(-\int_0^t f(s)\,\d J_s\right)=\exp\left(-\int_0^t \eta(f(s))\,\d s\right)\,.
\end{equation}
We observe that for $t\in[0,T_x)$ for the L\'evy process with secondary jump inputs no reflection has taken place.
As a consequence, for $t\in[0,T_x)$ we have $W^x_t=x+X_t$, and thus 
\begin{equation}
    W^x_t-W_t=(x+X_t)-(X_t+L_t)=x-L_t.
\end{equation}
Therefore, upon combining this with \eq{LT}, we have that
\begin{equation}\label{eq:diff}
A_x=\int_0^{T_x}h(W_t^x-W_t)\,\d t=\int_0^{T_x} h(x-L_t)\,\d t=\int_{(0,x]} h(x-y)\,T_{\d y}\,.
\end{equation}
Since $\{T_x\}_{x\geqslant 0}$ is a subordinator, the process $\{A_x\}_{x\geqslant 0}$ can be interpreted as a generalized L\'evy-driven shot-noise process.
Here the term `generalized' is used because in the standard form of shot noise \cite[Section 8.7]{Ross} the function $h(x)$ is the exponential function $e^{-r x}$ (for some decay rate $r>0$), and `L\'evy-driven' is used because in the standard form of shot noise this driving (increasing) Levy process is a compound Poisson process (cf.\ \cite{bkm19} for a general class of subordinator-driven shot-noise processes and induced infinite-server queues).

Since $\{T_y\}_{y\in[0,x]}$ is distributed like (the right continuous version of) $\{T_x-T_{x-y}\}_{y\in[0,x]}$, it also follows that, for any given $x>0$, the following distributional equality (in self-evident notation) applies:
\begin{equation}
A_x=\int_{(0,x]}h(x-y)\,T_{\d y}\stackrel{\rm d}{=}\int_{(0,x]} h(y)\,T_{\d y}\,.
\end{equation}
Recalling \eq{eta}, we thus arrive at the following result.

\begin{Theorem} \label{TH1} For any $\alpha\geqslant 0$, the LST of $A_x$ is given by
\begin{equation}
    {\mathbb E}\,e^{-\alpha A_x}=\exp\left({-\int_0^x \varphi^{-1}(\alpha h(y))\,{\rm d} y}\right)\,.
\end{equation}
\end{Theorem}

\begin{Remark}{\em 
In fact, with no additional effort we can also take $B_x:=\int_0^{T_x}g(W^x_t-W_t)\,\d t$ for some non-negative Borel function $g(\cdot)$, and observe that
for any $\alpha,\beta\geqslant 0$ the joint LST of $A_x$ and $B_x$ is given by
\begin{equation}
    {\mathbb E}\,e^{-\alpha A_x-\beta B_x}=\exp\left({-\int_0^x \varphi^{-1}(\alpha h(y)+\beta g(y))\,{\rm d} y}\right)\,.\label{R1}
\end{equation}
In particular, when plugging in $g(\cdot)\equiv 1$, we obtain the (joint) LST of $A_x$ and $T_x$. }$\hfill\Diamond$
\end{Remark}

As a consequence of \eqref{R1}, due to
\begin{equation}
(\varphi^{-1})'(0)=\frac{1}{\varphi'(0)}\quad\text{and}\quad (\varphi^{-1})''(0)=-\frac{\varphi''(0)}{(\varphi'(0))^3}
\end{equation}
and also recalling that $\varphi'(0)=-{\mathbb E}\,X_1$, $\varphi''(0)=\text{Var}\,X_1$, $(\varphi^{-1})'(0)={\mathbb E}\,T_1$ and $(\varphi^{-1})''(0)=-\text{Var}\,T_1$, 
it follows by differentiation (with respect to both $\alpha$ and $\beta$) and some standard algebra that
\begin{equation}
    {\mathbb E}\,A_x=\frac{1}{\varphi'(0)}\int_0^xh(y)\,\d y
\end{equation}
and that
\begin{equation}
    \text{Cov}(A_x,B_x)=\frac{\varphi''(0)}{(\varphi'(0))^3}\int_0^xh(y)g(y)\,\d y\,.
\end{equation}
Therefore, taking $g(\cdot)\equiv h(\cdot)$ gives
\begin{equation}
    \text{Var}\,A_x=\frac{\varphi''(0)}{(\varphi'(0))^3}\int_0^xh^2(y)\,\d y
\end{equation}
and with $g(\cdot)\equiv 1$ we have
\begin{equation}
    \text{Cov}(A_x,T_x)=\frac{\varphi''(0)}{(\varphi'(0))^3}\int_0^xh(y)\,\d y\,.
\end{equation}
The above expressions lead to the remarkable fact that the correlation coefficient between $A_x$ and $B_x$ is independent of the driving spectrally-positive L\'evy process $X$: it is given by
\begin{equation}
\text{Corr}(A_x,B_x)=\frac{\int_0^xh(y)g(y)\d y}{\sqrt{\int_0^xh^2(y)\,\d y\,\int_0^x g^2(y)\,\d y}}\,.
\end{equation}
As a consequence, the correlation coefficient between $A_x$ and $T_x$ reads
\begin{equation}
\text{Corr}(A_x,T_x)=\frac{\frac{1}{x}\int_0^xh(y)\d y}{\sqrt{\frac{1}{x}\int_0^xh^2(y)\d y}}\,.
\end{equation}

\begin{Remark}
   {\em 
When picking $h(x)=x$, we obtain results for the area between the processes $W^x$ and $W$ until time $T_x$:
\begin{align}
{\mathbb E}\,A_x=\frac{1}{\varphi'(0)}\frac{x^2}{2},\:\:\:\:
\text{Var}\,A_x=\frac{\varphi''(0)}{(\varphi'(0)^3}\frac{x^3}{3},\:\:\:\:
\text{Cov}(A_x,T_x)=\frac{\varphi''(0)}{(\varphi'(0))^3}\frac{x^2}{2},
\end{align}
so that
$\text{Corr}(A_x,T_x)=\frac{1}{2}{\sqrt{3}}$. 
Note that ${\mathbb E}\,A_x$ scales as $x^2$, as could have been anticipated. 
Evidently, by applying \eqref{R1} in principle any (joint) moment involving $A_x$ and $T_x$ can be expressed in terms of high-order derivatives of $\varphi(\cdot)$. 
It can be verified that when evaluating ${\mathbb E}\,[A_x^m\,T_x^n]$ one obtains a monomial of degree $2m+n$.
}\hfill$\Diamond$
\end{Remark}

\begin{Remark}\rm
Note that when $\varphi''(0)<\infty$, both ${\mathbb E}\,T_x$ and $\text{Var}\,T_x$ scale linearly in $x$:
\begin{equation}
{\mathbb E}\,T_x=\frac{x}{|{\mathbb E}\,X_1|}\quad \text{and}\quad \text{Var}\,T_x=\frac{\text{Var}\,X_1}{|{\mathbb E}\,X_1|^3}\, x\,.
\end{equation}
Recall that, with $\sigma(Y):=\sqrt{\text{Var}\,Y}$ for some random variable $Y$,
\begin{equation}
    \frac{T_x-{\mathbb E}\,T_x}{\sigma(T_x)}=\left.\left({T_x-\frac{x}{|{\mathbb E}\,X_1|}}\right)\right/\left({\sqrt{\frac{\text{Var}\,X_1}{|{\mathbb E}\,X_1|^3}\,x}}\right)
\end{equation}
converges in distribution to a standard normal random variable as $x\to\infty$. Although this particular case is a very simple exercise (either via the central limit theorem or via convergence of corresponding exponents), for a thorough study of central limit theorems for L\'evy processes, see \cite{dm01}.

This reminds us that for any renewal counting process $N\equiv \{N_t\}_{t\geqslant 0}$ we have an identical limit, in the sense that when the increments of the renewal process are distributed like some non-negative (but not a.s.\ zero) random variable $Y$ and in addition ${\mathbb E}\,Y^2<\infty$, then (e.g., \cite[Ch.~V, Prop. 6.3]{asmussen})
\begin{equation}
    \left.\left({N_t-\frac{t}{{\mathbb E}\,Y}}\right)\right/
    \left({\sqrt{\frac{\text{Var}\,Y}{({\mathbb E}\,Y)^3}\, t}}\right)
\end{equation}
also converges in distribution to a standard normal random variable as $t\to\infty$.
Notice though that there is a crucial distinction between both results: the process $T$ has stationary independent increments, whereas $N$ typically does not (unless it is a Poisson process).
\hfill$\Diamond$
\end{Remark}

\begin{Remark} \rm
    We note that identity \eq{diff} is valid whenever $X$ has no negative jumps and $X_t\to-\infty$ a.s.\ as $t\to\infty$. 
    That is, for this property to hold we do not need to assume that $X$ is a L\'evy process. 
    We also note that because of that, whenever the assumptions are such that $T_x$ has stationary (but not necessarily independent) increments with ${\mathbb E}\,T_1<\infty$, then $ET_x=ET_1\, x$, so that
    \begin{equation}
       {\mathbb E}\,A_x={\mathbb E}\,T_1\int_0^xh(y)\,{\rm d}y, 
    \end{equation} just like in the case with independent increments. For other moments such a property, of course, does not hold. 
\hfill$\Diamond$
\end{Remark}

\begin{Remark}\label{rem:longrun}\rm
    By applying results from regenerative theory (e.g., \cite[Ch. VI, Th.~3.1]{asmussen}), we know that
    \begin{equation}
\lim_{t\to\infty}\frac{1}{t}\int_0^t h(W_t^x-W_t)\,\d t=\frac{{\mathbb E}\,A_x}{{\mathbb E}\,T_x}=\frac{1}{x}\int_0^xh(y)\,\d y={\mathbb E}\,h(xU)\,,
\end{equation}
where $U\sim\text{Uniform}(0,1)$ (so that $xU\sim \text{Uniform}(0,x)$).
It may seem surprising that this long-run average limit is independent of the particular driving L\'evy process $X$ and depends only on the function $h(\cdot)$ and secondary jump size $x$. 
However, a similar phenomenon can be seen from the fact \cite{kw92} that the steady-state distribution of $W^x$ is a convolution of the steady state distribution of $W$ and that of $xU$, where the latter is independent of the L\'evy process $X$. 
Nevertheless, the proofs of these two results are different and, although it makes sense intuitively, it is not obvious how one of these two phenomena can be easily concluded from the other. 
In the special case that $h(x)=x$, both results will imply that the difference between the means of the steady-state distributions of $W^x$ and $W$ is ${\mathbb E}\,[xU]=x/2$ (provided that $\int_{(1,\infty)}x^2\,\nu(\d x)<\infty$, otherwise both means are infinite).
\hfill$\Diamond$
\end{Remark}

\begin{Remark}\rm
So far we have considered the area between $W^x$ and $W$, but the reasoning extends to the 
area of a function of the difference of $W^y_t$ and $W^x_t1_{[0,T_x)}(t)+W_t1_{[T_x,T_y)}(t)$ until time $T_y$, for $0\leqslant x\leqslant y$.
That is, the area of a function of the difference between $y+X_t$ and the reflection of $x+X_t$ until time $T_y$.
Note that this area is a sum of two independent random variables: the first, corresponding to the time interval until $T_x$, is equal to $h(y-x)T_x$, whereas the second, corresponding to the time interval between $T_x$ and $T_y$, is distributed like $A_{y-x}$. 
This claim is a direct consequence of the strong Markov property of (right continuous) L\'evy processes.

We proceed by evaluating the LST\,s of the two components constituting
\begin{equation}A_{x,y} :=\int_0^{T_y} h(W_t^y-(W^x_t1_{[0,T_x)}(t)+W_t1_{[T_x,T_y)}(t)))\,{\rm d}t\stackrel{\rm d}{=} h(y-x)T_x + A_{y-x},\end{equation}
with the two random quantities in the right-most expression being independent.
The LST of $h(y-x)T_x$ is 
\begin{equation}
{\mathbb E}\,\exp\left({-\alpha h(y-x)\,T_x}\right)=\exp\left({-\varphi^{-1}(\alpha h(y-x))\,x}\right),
\end{equation}
and the mean and variance are 
\begin{equation}\frac{h(y-x)\,x}{|{\mathbb E}\,X_1|}\ \text{and}\ 
\frac{h^2(y-x)\,x\,\text{Var}\,X_1}{|{\mathbb E}\,X_1|^3}\,,
\end{equation} 
respectively. The LST of $A_{y-x}$ follows from Theorem \ref{TH1}.
Upon combining the above, we conclude that the LST of $A_{x,y}$ is given by, for any $\alpha\geqslant 0$,
\begin{equation}
{\mathbb E}\,e^{-\alpha A_{x,y}}=
\exp\left(-\varphi^{-1}(\alpha h(y-x))\,x+\int_0^{y-x}\varphi^{-1}(\alpha h(z))\,{\rm d}z\right)
\end{equation}
Clearly,
\begin{align}
    {\mathbb E}A_{x,y}&=h(y-x){\mathbb E} T_x+{\mathbb E} A_{y-x}=\frac{h(y-x)x+\int_0^{y-x}h(s)\text{d}s}{|{\mathbb E} X_1|}\nonumber\\
    \text{Var}\,A_{x,y}&=h^2(y-x)\text{Var}\,T_x+\text{Var}\,A_{y-x}\\
    &=\frac{\text{Var}\, X_1}{|{\mathbb E} X_1|^3}\left(h^2(y-x)x+\int_0^{y-x} h^2(s)\text{d}s\right)\,,\nonumber
\end{align}
for $0\leqslant x\leqslant y$. \hfill$\Diamond$
\end{Remark}

\begin{Remark}\em
It is straightforward to find the joint LST for the finite-dimensional distribution of $A_x$. For $0=s_0<s_1<\ldots<s_n$ and $x_i=s_i-s_{i-1}$ for $i=1,\ldots,n$ we can write
\begin{equation}
    A_{s_i}=\sum_{j=1}^i\int_{(0,x_j]}h(s_i-s_{j-1}-t)T_{s_{j-1}+\d t}
\end{equation}
and thus
\begin{align}
    \sum_{i=1}^n\alpha_iA_{s_i}&=\sum_{i=1}^n\alpha_i\sum_{j=1}^i\int_{(0,x_j]}h(s_j-s_{j-1}-t)T_{s_{j-1}+\d t}\nonumber \\ &=\sum_{j=1}^n\int_{(0,x_j]}\left(\sum_{i=j}^n\alpha_i h(s_i-s_{j-1}-t)\right) T_{s_{j-1}+\d t}
\end{align}
where we note that $T$ being a subordinator implies that the summands are independent and distributed as
\begin{equation}
    \int_{(0,x_j]}\left(\sum_{i=j}^n\alpha_i h(s_i-s_{j-1}-t)\right)T_{\d t}\,.
\end{equation}
Therefore, for any non-negative $\alpha_1,\ldots,\alpha_n$, the joint LST is given by
\begin{align}
    {\mathbb E}\exp\left(-\sum_{i=1}^n\alpha_i A_{s_i}\right)&=\prod_{j=1}^n\exp\left(-\int_0^{x_j}\varphi^{-1}\left(\sum_{i=j}^n\alpha_ih(s_i-s_{j-1}-t)\right)\d t\right)\\
    &=\prod_{j=1}^n\exp\left(-\int_0^{x_j}\varphi^{-1}\left(\sum_{i=j}^n\alpha_ih(s_i-s_j+t)\right)\d t\right)\,.\nonumber
\end{align}
Along the same lines it follows that
\begin{align}
   {\mathbb E}\exp\left(-\sum_{i=1}^n\alpha_i A_{s_i}-\sum_{i=1}^n\beta_iT_{s_i}\right)
   &=\prod_{j=1}^n\exp\left(-\int_0^{x_j}\varphi^{-1}\left(\sum_{i=j}^n\left(\alpha_ih(s_i-s_j+t)+\beta_i\right)\right)\d t\right)\,,
\end{align}
for non-negative  $\alpha_i,\beta_i$, $i=1,\ldots,n$; cf.\ \eqref{R1}. \hfill$\Diamond$
\end{Remark}

\begin{Remark}\rm 
    One may wonder whether it would be easier to study the stochastic properties of the areas $\int_0^{T_x}W^x_t\,\d t$ and $\int_0^{T_x}W_t\,\d t$ separately and then compare, rather than considering their difference. 
    As mentioned in the introduction,  both turn out to be intrinsically more difficult problems, which have so far not been successfully handled in the generality of a spectrally-positive L\'evy process (but see the recent results on moments of $\int_0^{T_x}W^x_t\,\d t$ in \cite{GJM} for a subclass of spectrally-positive L\'evy processes). 
    It is not hard to show that $\int_0^{T_x}W_s\,\d s$ (indexed by $x$) is a subordinator, but the problem of identifying its exponent seems to be quite hard in general. \hfill$\Diamond$
\end{Remark}

\section{Higher moments}\label{S4}
Suppose that $f(\alpha)=e^{-g(\alpha)}$ is an LST of some non-negative random variable $Y$.
This holds in the case that $g(\cdot)$ is (the negative of) a cumulant moment generating function, i.e., $g(\alpha) =-\log{\mathbb E} e^{-\alpha Y}$.
Then both $f(\cdot)$ and $g(\cdot)$ are differentiable of all orders, $f(0)=1$, $g(0)=0$ and we have that
\begin{equation}
f'(\alpha)=-g'(\alpha)f(\alpha)\,.
\end{equation}
Thus,
\begin{equation}
(-1)^{n+1}f^{(n+1)}(\alpha)=\sum_{k=0}^n \binom{n}{k}(-1)^k g^{(k+1)}(\alpha)(-1)^{n-k}f^{(n-k)}(\alpha)\,.
\end{equation}
Denoting $\mu_0:=1$, $c_0:=0$ and (for $k=1,2,\ldots$), $\mu_k:=(-1)^kf^{(k)}(0)$ and $c_k:=(-1)^{k-1}g^{(k)}(0)$, we immediately see that, whenever $g'$ is completely monotone (as in our case), $\mu_n$ can be found recursively:
\begin{equation}\label{eq:moments}
\mu_{n+1}=\sum_{k=0}^n \binom{n}{k}c_{k+1}\mu_{n-k}
\end{equation}
which, by monotone convergence holds regardless of the finiteness of $c_k$ and $\mu_k$ appearing in the right-hand side. 
It follows by induction that $\mu_n<\infty$ if and only if $c_n<\infty$, and that both imply that $\mu_i,c_i<\infty$ for $i\in\{0,\ldots, n\}$.

In view of Theorem \ref{TH1}, our interest is in the special case where $g(\alpha)=\int_0^x\varphi^{-1}(\alpha h(y))\,\d y$. In this case we have that
\begin{equation}
c_k=(-1)^{k-1}(\varphi^{-1})^{(k)}(0)\int_0^xh^k(y)\,\d y\,.
\end{equation}
In particular, when $h(x)=x$ we have that $\int_0^xh^k(y)\,\d y=x^{k+1}/(k+1)$.
As was mentioned earlier, the derivatives of $\varphi^{-1}(\cdot)$ at zero may be inferred from the derivatives of $\varphi(\cdot)$ at zero via the Fa\`a di Bruno formula.

\begin{Remark}\rm
Moments can also be computed when ordering a random amount of secondary work. Namely, consider a general $h(\cdot)$, and let $x$ be replaced by a non-negative, finite mean random variable $\xi$ with ${\mathbb P}(\xi=0)<1$, which is independent of everything else.
Then
\begin{equation}
{\mathbb E}\int_0^\xi h^k(y)\,\d y={\mathbb E}\xi\,{\mathbb E}h^k(\xi_e)\,,
\end{equation}
where $\xi_e$ has the stationary excess lifetime distribution associated with $\xi$ (with density $f_{\xi_e}(t)={\mathbb P}(\xi>t)/{\mathbb E}\xi$). 
When $h(x)=x$, we thus find that
\begin{equation}
{\mathbb E}\xi\,{\mathbb E}\xi_e^k=\frac{{\mathbb E}\xi^{k+1}}{k+1}\,.
\end{equation}
When in addition $\xi\sim\exp(\lambda)$,
\begin{equation}
{\mathbb E}\xi\,{\mathbb E}\xi_e^k=\lambda^{-1} {\mathbb E}\xi^k=\frac{{\mathbb E}\xi^{k+1}}{k+1}=\frac{k!}{\lambda^{k+1}}\,,
\end{equation}
in which case \eq{moments} implies that, with $\tilde \mu_n=\mathbb{E}A^n_\xi$ and $\tilde c_n=(-1)^{n-1}(\varphi^{-1})^{(n)}(0)$,
\begin{equation}
\tilde\mu_{n+1}=n!\sum_{k=0}^n\frac{(k+1)\tilde c_{k+1}}{\lambda^{k+2}}\frac{\tilde\mu_{n-k}}{(n-k)!}\,,
\end{equation}
for $n\in{\mathbb N}$.\hfill$\Diamond$
\end{Remark}

\section{Functional Gaussian limit theorem for $A_x$}\label{S5}
The scope of this section is a functional Gaussian limit theorem for a time-scaled version of the process $A\equiv \{A_{x}\}_{x\geqslant 0}$.

Recall that $h(\cdot)$ is said to be a {\it regularly varying function} with index $\alpha\geqslant 0$ if it has the form $x^\alpha \ell(x)$, where $\ell(\cdot)$ is slowly varying (i.e., $\ell(ax)/\ell(x)\to 1$ as $x\to\infty$ for all $a>0$). The following result has the flavor of a functional central limit theorem, in that there is convergence to a Gaussian limit, but the normalization differs from the usual $\sqrt{n}$.

\begin{Theorem}\label{th2}
Let $h:[0,\infty)\to [0,\infty)$ be a right-continuous, nondecreasing function that is regularly varying (at infinity) with some index $\alpha\geqslant 0$. If $\mbox{\rm Var}\,X_1=\varphi''(0)<\infty$ and $-{\mathbb E}X_1=\varphi'(0)>0$,
\begin{equation}\label{eq:hn}
\frac{1}{h(n)n^{1/2}}\left(A_{nx}-\frac{1}{\varphi'(0)}\int_0^{nx}h(s)\, {\rm d} s\right)   
\end{equation}
converges weakly in $(D[0,\infty),J_1)$ (with $J_1$ referring to the {\it  Skorokhod $J_1$ topology}) to
\begin{equation}
A^*_x\equiv \int_0^x(x-s)^\alpha \,{\rm d} B_s=\alpha\int_0^x (x-s)^{\alpha-1}\,B_s\,{\rm d} s,
\end{equation}
where $B$ is a driftless Brownian motion with variance coefficient ${\rm Var}\,T_1=\varphi''(0)/(\varphi'(0))^3$. 
\end{Theorem}

\begin{Remark}
When $h:[0,\infty)\to [0,\infty)$ is right continuous, nondecreasing, with $\ell(x)=h(x)/x^\alpha\to\ell\in(0,\infty)$ as $x\to\infty$, then it is also regularly varying and $h(n)n^{1/2}$ in \eq{hn} may be replaced by $\ell n^{\alpha+1/2}$.
\end{Remark}

\noindent{\it Proof.}
In \cite[Theorem~1.1]{ir20} a functional limit theorem is presented for shot-noise processes of the form $\int_0^th(t-s)\,\d N_s$, where $N$ is a counting process with respect to some point process. 
Upon checking the underlying argumentation, it turned out that the counting process assumption is not necessary, which was confirmed by the authors \cite{ir23}.
This observation has the direct implication that \cite[Theorem~1.1]{ir20} extends to the setting we are considering now. In particular, one may infer from their results that a sufficient condition for Theorem~\ref{th2} to hold is that
\begin{equation}\label{eq:wct}
\frac{T_{nx}-n\,{\mathbb E}T_x}{\sqrt{n}}=\frac{T_{nx}-nx/\varphi'(0)}{\sqrt{n}}
\end{equation}
converges (as $n\to\infty$) weakly in $(D[0,\infty),J_1)$ to a driftless Brownian motion with variance coefficient $\text{Var}\,T_1$.

Clearly, for $x=1$, the right and left hand sides of \eq{wct} converge in distribution to a normal distribution with mean zero and variance $\text{Var}\, T_1$. By \cite[Theorem~1]{p07} (see also \cite[Th.~2.7]{skorokhod}) this implies the functional limit theorem for \eq{wct} and therefore, recalling that
\begin{equation}{\mathbb E}A_x=\frac{1}{\varphi'(0)}\int_0^xh(s)\,\d s
\end{equation}
for any non-negative  Borel function $h(\cdot)$, we are done.\hfill $\Box$

\medskip

It should be mentioned that for each fixed $x$, $\int_0^x(x-s)^\alpha \d B_s$ is distributed like $\int_0^x s^\alpha \d B_s$, which is normally distributed with expected value zero and variance 
\begin{equation}\text{Var}\,T_1\,\int_0^x s^{2\alpha}\,\d s=\text{Var}\,T_1\,\frac{x^{2\alpha+1}}{2\alpha+1}\,.
\end{equation}
It can also be checked that, for $x,y\ge 0$, the autocovariance for the Gaussian process $\{A^*_x\}_{x\geqslant 0}$ is
\begin{equation}
    \text{Cov} (A^*_x,A^*_{x+y})=\text{Var}\,T_1\, \int_0^x[(x-s)(x+y-s)]^\alpha \d s=\text{Var}\,T_1 \int_0^x [s(y+s)]^\alpha \d s\,.
\end{equation}
Hence, as expected, the correlation coefficient, depends only on $x,y,\alpha$. 

\section{An inventory application}\label{S6}
Here we assume that we would like to keep the server constantly busy by having it perform a secondary job when it is not processing the primary material.
As mentioned earlier, $W^x_t-W_t$ models this secondary (lower priority) content in the system at time $t$.
The holding cost for this secondary content that we pay at time $t$ is $h(W^x_t-W_t)\text{d}\,t$ where $h$ is non-decreasing and right-continuous. Then, as was observed in Remark \ref{rem:longrun},  the long-run average extra cost for ordering $x$ at the beginning of each cycle is
\begin{equation}
\frac{1}{x}\int_0^xh(s)\,\d s\,.
\end{equation}
Let us assume that for each cycle we pay an additional setup cost $K>0$ for the shipment of $x$ secondary job units, which is assumed independent of the shipment size. Then the total long run average setup cost is $K/{\mathbb E}T_x=\varphi'(0)K/x$, and thus the long-run average total cost is
\begin{equation}\label{eq:lac}
\frac{\varphi'(0)K+\int_0^xh(s)\,\d s}{x}\,.
\end{equation}
Denote $K':=\varphi'(0)K$ and consider the non-decreasing function
\begin{align}\label{eq:g}
g(x)&=\int_0^x(h(x)-h(s))\,\d s=\int_0^x\int_{(s,x]}h(\d t)\,\d s=\int_{(0,x]}\int_0^t\d s\,h(\d t)=\int_{(0,x]}t\,h(\d t)\,.
\end{align}
Assume that $\{x\,|\,g(x)\geqslant K'\}$ is nonempty and let
\begin{equation}\label{eq:xstar}
x^*:=\inf\{x\,|\,g(x)\geqslant K'\}\,.
\end{equation}
Then, by results for generalized inverses, $g(x)\geqslant K'$ (resp., $g(x)<K'$) if and only if $x\geqslant x^*$ (resp., $x<x^*$).
Now, for every $x<x^*$ we have that $g(x)<K'$. Recalling that $g(x)=xh(x)-\int_0^xh(s)\,\d s$, this implies that for every $0\leqslant s\leqslant x$ we have
\begin{equation}
h(s)\leqslant h(x)<\frac{K'+\int_0^xh(s)\d s}{x}\,.
\end{equation}
This implies that for $y<x$
\begin{equation}
\frac{1}{x-y}\int_y^xh(s)\d s\leqslant h(x)<\frac{K'+\int_0^xh(s)\d s}{x}\,.
\end{equation}
It is a straightforward exercise to show that this is equivalent to, for $y<x$,
\begin{equation}
\frac{K'+\int_0^xh(s)\d s}{x}<\frac{K'+\int_0^yh(s)\d s}{y},
\end{equation}
entailing that the cost function is decreasing on $(0,x^*)$; similarly, it is non-decreasing on $(x^*,\infty)$. Noting that the function is continuous on $(0,\infty)$ and in particular at $x^*$, it follows that $x^*$ is a minimum.

When $\{x\,|\,g(x)\geqslant K'\}$ is empty, the cost function is strictly decreasing and therefore there is no finite minimal point. When $h(\cdot)$ is not a constant and has a density $h'(\cdot)$ for which $t\mapsto th'(t)$ is non-decreasing and unbounded (e.g.\ when $h(\cdot)$ is linear or convex), then this is not possible.  In particular when $h(t)=t^\gamma$ where $\gamma>0$ (even when $\gamma\in(0,1)$, in which case $h(\cdot)$ is strictly concave).

A related question is: what should be the per-unit cost of preventing the process for becoming negative (e.g., ordering from a more expensive source in order to satisfy demand) and the per-unit reward for secondary output be in order for us to prefer ordering $x^*$ every time the system become empty. If the penalty per unit time is $p$, then the long run average cost for the reflected case would be the (a.s.) limit of $pL(t)/t$ where $L(t)=-\inf_{0\leqslant s\leqslant t}X_s$, which equals $p\varphi'(0)$. If the reward per unit of secondary output is $r$, then since $rx$ is the reward per (regenerative) cycle, then the reward per unit time is $\frac{rx}{ET_x}=\varphi'(0)r$ for every $x>0$ and in particular for $x=x^*$.

This implies that
\begin{equation}
p^*=\frac{K+\frac{1}{\varphi'(0)}\int_0^{x^*}h(t)\,\d t}{x^*}-r
\end{equation}
is the break-even penalty in the sense that for $p<p^*$ we would prefer not to perform secondary work and for $p>p^*$ we would order up to $x^*$ every time the workload is exhausted.

When $h(t)=ct$, we have that $g(x)=\int_{(0,x]}th(\d t)=c\int_0^xt\d t=cx^2/2$ and thus we obtain the familiar expression 
\begin{equation}
    x^*=\sqrt{\frac{2\varphi'(0)K}{c}}\,.
\end{equation} In this case
\begin{equation}
\frac{\varphi'(0)K+\int_0^{x^*}h(t)\,\d t}{x^*}=\sqrt{2\varphi'(0)Kc}\,,
\end{equation}
which is also a solution of the textbook economic order quantity problem with setup cost $K$, demand rate $d=1/\varphi'(0)$ and linear holding cost $c$.

\medskip

We recall that the results of the present section is for a system in which the L\'evy process $X$ is interpreted as primary input, and the `spikes' of size $x$ as secondary input, where the primary input have priority over the secondary input. Therefore the area of a function of the difference between the two workload graphs reflects the performance perceived by the secondary input.

This idea can be further extended, as follows. Define, for $j\in\{1,\ldots,m\}$ for some $m\in{\mathbb N}$, $s_j:=\sum_{i=1}^jx_i$, where $x_i\ge 0$. Now consider the model in which the spike has size $s_m$. This could model a situation with $m+1$ priority classes, with $W$ representing the workload of the highest priority class; order them such that class $i$ has priority over class $j$ if $i<j$. Whenever the systems becomes idle, a new spike of size $s_m$ enters. 
For $k\in\{ 1,\ldots, m-1\}$, the content corresponding to the $(k+1)$-st priority class during $[0,T_{s_{k-1}})$ is $x_k$. During $[T_{s_{k-1}},T_{s_k})$ it is $W_{t-T_{s_{k-1}}}-W^{x_k}_{t-T_{s_{k-1}}}$, and during $[T_{s_k},T_{s_m})$ it is zero.
In this situation we can attempt to perform the same optimization as before, but now over the vector $(x_1,\ldots,x_m)\in{\mathbb R}_+^m$, so as to determine the optimal $x_1^\star,\ldots,x_m^\star$.

In general this multivariate optimization is not so easy (at least not analytically), other than for the linear case, which will be explained shortly. If $h:\mathbb{R}^m\to\mathbb{R}$ is a nondecreasing holding cost function, then a similar derivation that leads to \eq{lac}, implies that the long run average cost associated with the secondary priorities has the form
\begin{equation}\label{eq:mtac}
\frac{\varphi'(0)K+\sum_{i=1}^m\int_0^{x_i}h(0,\ldots,0,s,x_{i+1},\ldots,x_m)\,\text{d}s}{s_m}\,.
\end{equation}
If we insist that the total order is $x$ and that the proportion ordered from each class $i$ is some fixed $p_i$, then \eq{mtac} may be written as
\begin{equation}
\frac{\varphi'(0)K+\int_0^x \sum_{i=1}^m p_ih(0,\ldots,0,p_is,p_{i+1}x,\ldots,p_mx)\,\text{d}s}{x}\,.
\end{equation}
This is not the same problem as the one we solved earlier, since now there is also $x$ in the integrand. However, if the assumptions are such that the function
\begin{align}
H(x)&\equiv\int_0^x \sum_{i=1}^m p_ih(0,\ldots,0,p_is,p_{i+1}x,\ldots,p_mx)\,\text{d}s\nonumber\\
&=x\sum_{i=1}^m \int_0^{p_i}h(0,\ldots,0,sx,p_{i+1}x,\ldots,p_mx)\,\text{d}s
\end{align}
is convex in $x$ on $(0,\infty)$, then the optimal $x$ has the same structure as in \eq{g} and \eq{xstar} where $h$ is replaced by the right continuous version of the density of $H$.

When $h(x_1,\ldots,x_m)=\sum_{i=1}^mh_i(x_i)$, then it is easy to check that \eq{mtac} becomes
\begin{equation}
\frac{\varphi'(0)K+\sum_{i=1}^m\left(\int_0^{x_i}h_i(s)\,\text{d}s+s_{i-1}h_i(x_i)\right)}{s_m}\,,
\end{equation}
where $s_0=0$. With $x_i=p_ix$ (so that $s_m=x$) and $F_j=\sum_{i=1}^j p_i$ (with $F_0=0$), this becomes
\begin{equation}\label{eq:sepp}
\frac{\varphi'(0)K}{x}+\sum_{i=1}^m\left(\int_0^{p_i}h_i(xs)\,\text{d}s+F_{i-1}\,h_i(p_ix)\right)\,.
\end{equation}
Therefore, it is possible that under some assumptions such a problem may be solved in two stages. First fix $x$ and find the optimal minimizing proportions. If the expression that is obtained is some well behaved function of $x$, then  we may optimize with respect to $x$ either analytically or, if that is not possible, numerically. 

One case that gives an explicit (albeit boring) solution is when $h_i(t)=c_it$ for $c_i>0$ for all $i=1\ldots,m$. In this case it is straightforward to check that \eq{sepp} becomes

\begin{equation}
\frac{\varphi'(0)K}{x}+\frac{x}{2}\sum_{i=1}^m(F_i^2-F_{i-1}^2)\,c_i\ge \frac{\varphi'(0)K}{x}+\frac{xc}{2}  \,,
\end{equation}
where $c\equiv\min_{1\leqslant k\leqslant m}c_k$.  It therefore follows that optimal proportions (not necessarily unique) are such that $p_i=0$ for each $i$ for which $c_i>c$ and $x=\sqrt{2\varphi'(0)K/c}$. This is the same as writing that $x_i=0$ for each $i$ for which $c_i>c$ and that $\sum_{i|c_i=c}x_i=\sqrt{2\varphi'(0)K/c}$.


\begin{thebibliography}{99}

\bibitem{applebaum} Applebaum, D. (2009). {\em L\'evy Processes and Stochastic Calculus.} Cambridge.

\bibitem{asmussen} Asmussen, S. (2003). {\em Applied probability and queues, 2nd Edition.} Springer.

\bibitem{bertoin} Bertoin, J. (1996). {\em L\'evy Processes.} Cambridge.


\bibitem{Borovkov2003}
Borovkov, A. A., Boxma, O. J.,  Palmowski, Z. (2003). On the integral of the workload process of the single server queue. \textit{Journal of Applied Probability}, \textbf{40}, 200-225.

\bibitem{bk14} Boxma O., Kella O. (2014). Decomposition results for stochastic storage processes and queues with alternating Lévy inputs. {\em 
Queueing Systems,} {\bf 77,} 97-112.

\bibitem{bkm19} Boxma O., Kella O., Mandjes M. (2019). Infinite-server systems with Coxian arrivals. {\em Queueing Systems} {\bf 92,} 233-255.

\bibitem{dm01} Doney, R. A. and Maller, R. A. (2002). Stability and attraction to normality for L\'evy processes at zero and at infinity. {\em Journal of Theoretical Probability,} {\bf 15,} 751-792.

\bibitem{DM}
Debicki, K., Mandjes, M. (2015). \textit{Queues and L\'evy fluctuation theory}. Springer.

\bibitem{GJM} Glynn, P., Jacobovic, R., Mandjes, M. (2023). 
Moments of polynomial functionals in L\'evy-driven queues with secondary jumps. {\it Submitted.}


\bibitem{Iglehart1971}
Iglehart, D.  (1971). Functional limit theorems for the queue GI/G/1 in light traffic. \textit{Advances in Applied Probability}, \textbf{3}, 269-281.


\bibitem{ir20} Iksanov, A. and Rashytov, B. (2020). A functional limit theorem for general shot noise processes. {\em Journal of Applied Probability} {\bf 57,} 280-294.


\bibitem{ir23} Iksanov, A. and Rashytov, B. Private communication.

\bibitem{jv83} Jurek, Z. J. and Vervaat, W. (1983). An integral representation for sefdecomposable Banach space valued random variables. {\em Zeitschrift f\"ur Wahrsheinlichkeitstheorie und Verwandte Gebiete.} {\bf 62,} 247-262.

\bibitem{k06} Kella, O. (2006). Reflecting thoughts. {\em Statistics and Probability Letters} {\bf 76,} 1808-1811.

\bibitem{kw92} Kella O., Whitt W. (1992). Useful martingales for stochastic storage processes with Lévy input. {\em Journal of Applied Probability,} {\bf 29,} 396-403.

\bibitem{kw91} Kella O., Whitt W. (1991). Queues with server vacations and Lévy processes with secondary jump input. {\em The Annals of Applied Probability,} {\bf 1,} 104-117.

\bibitem{kw90} Kella O., Whitt W. (1990). Diffusion approximations for queues with server vacations. {\em Advances in Applied Probability,} {\bf 22,} 706-729.

\bibitem{kyprianou} Kyprianou, A. E. (2006). {\em Introductory Lectures on Fluctuations of L\'evy Processes with Applications.} Springer.

\bibitem{lukacs} Lukacs, E. (1975). {\em Stochastic convergence, 2nd edition.} Academic Press.

\bibitem{p07} Permiakova, E. (2007). Functional limit theorems for L\'evy processes and their almost-sure versions. {\em Lithuanian Mathematical Journal} {\bf 47,} 67-78.

\bibitem{Ross}
Ross, S. (1996). {\em Stochastic Processes}. Wiley.

\bibitem{sato} Sato, K. (1999). {\em L\'evy Processes and Infinitely Divisible Distributions.} Cambridge. 

\bibitem{skorokhod} Skorokhod, A. B. (1957). Limit theorems for stochastic processes with independent increments. {\em Theory of Probability and its Applications.} {\bf 2,} 138-171.

\end{thebibliography}
\end{document}